\documentclass[12pt]{article}

\title{On improved estimation in a conditionally Gaussian regression}

\author {Pchelintsev Evgeny\thanks{
Department of Mathematics and Mechanics, Tomsk State University, Lenin str. 36, 634050 Tomsk, Russia, and Laboratoire de
Math\'ematiques Rapha\"{e}l Salem, UMR 6085 CNRS, Universit\'e de Rouen, Avenue de l'Universit\'e BP.12, 76800 Saint Etienne du Rouvray Cedex, France, e-mail: evgen-pch@yandex.ru}}

\usepackage{amssymb}
\usepackage{amsfonts}
\usepackage{amsmath}
\usepackage{amsthm}
\usepackage{enumerate}
\usepackage{multicol}
\usepackage{graphicx}

\newtheorem{theorem}{Theorem}[section]
\newtheorem{proposition}[theorem]{Proposition}
\newtheorem{lemma}[theorem]{Lemma}

\newtheorem{remark}{Remark}[section]
\newtheorem{corollary}[theorem]{Corollary}

\def\text#1{\hbox{#1}}
\def\proof{{\noindent \bf Proof. }}
\def\endproof{\hspace*{130mm}\mbox{\ $\qed$}}

\def\build #1_#2{\mathrel{\mathop{\kern 0pt #1}\limits_{#2}}}
\newcommand{\zs}[1]{{\mathchoice{#1}{#1}{\lower.25ex\hbox{$\scriptstyle#1$}}
{\lower0.25ex\hbox{$\scriptscriptstyle#1$}}}}

\begin{document}

\maketitle

\begin{abstract}
The paper considers the problem of estimating a $p\geq2$\ dimensional mean vector of a multivariate conditionally normal distribution under quadratic loss. The problem of this type arises when estimating the parameters in a continuous time regression model with a non-Gaussian Ornstein--Uhlenbeck process driven by the mixture of a Brownian motion and a compound Poisson process. We propose a modification of the James--Stein procedure of the form $\theta^*(Y)=(1-c/\|Y\|)Y,$\ where $Y$\ is an observation and $c>0$\ is a special constant. This estimate allows one to derive an explicit upper bound for the quadratic risk and has a significantly smaller risk than the usual maximum likelihood estimator for the dimensions $p\geq2$.\ This procedure is applied to the problem of parametric estimation in a continuous time conditionally Gaussian regression model and to that of estimating the mean vector of a multivariate normal distribution when the covariance matrix is unknown and depends on some nuisance parameters.

\end{abstract}

\vspace*{5mm} \noindent {\sl Keywords}: Conditionally Gaussian regression model; Improved estimation; James--Stein procedure; Non-Gaussian Ornstein--Uhlenbeck process.

\vspace*{5mm} \noindent {\sl AMS 1991 subject classifications}:
Primary: 62C20; Secondary: 62C15

\bibliographystyle{plain}
\renewcommand{\columnseprule}{.1pt}

\section{Introduction}\label{sec:int}

In 1961, James and Stein \cite{JS} considering the problem of estimating the mean vector $\theta$\ of a $p$-dimensional\ normal distributed random vector $Y$\ with a covariance matrix $I_p$\ introduced an estimator
\begin{equation}\label{sec:int.1}
\hat{\theta}_{JS}=\left(1-\frac{p-2}{\|Y\|^2}\right)Y
\end{equation}
which outperforms the maximum likelihood estimate (MLE)
\begin{equation}\label{sec:int.2}
\hat{\theta}_{ML}=Y
\end{equation}
for dimension $p\geq3$,\ under the common quadratic risk 
\begin{equation}\label{sec:int.3}
R(\theta,\hat{\theta})=\mathbf{E}_{\theta}\|\theta-\hat{\theta}\|^{2},
\end{equation}
in the sense that for all parameter values $\theta$\ 
\begin{equation*}
R(\theta,\hat{\theta}_{JS})<R(\theta,\hat{\theta}_{ML}).
\end{equation*}
This unexpected result draw a great interest of mathematical statisticians and stimulated a number of authors to contribute to the theory of improved estimation by extending the problem of James and Stein in different directions to more general models with unknown covariance matrix and considering other types of estimates (see \cite{BeHa,EfMo,Gl2,St2} for more details and other references). A considerable effort has been directed towards the problems of improved estimation in non-Gaussian models with the spherically symmetric distributions (see \cite{Fou,FoSt2}) and in the non-parametric regression models \cite{FoPe}.

Now the James--Stein estimator and other improved shrinkage estimators are widely used in econometrics and the problems associated with the signal processing.

In this paper we will consider the problem of estimating the mean in a conditionally Gaussian distribution. Suppose that the observation $Y$\ is a $p$-dimensional random vector which obeys the equation
\begin{equation}\label{sec:int.4}
Y=\theta+\xi,
\end{equation}
where $\theta$\ is a constant vector parameter, $\xi$\ is a conditionally Gaussian random vector with a zero mean and the covariance matrix $\mathcal{D}(\mathcal{G})$,\ i.e. $Law(\xi|\mathcal{G})=\mathcal{N}_{p}(0,\mathcal{D}(\mathcal{G}))$,\ where $\mathcal{G}$\ is some fixed $\sigma$-algebra.\

We propose to consider a shrinkage estimator of the form
\begin{equation}\label{sec:int.5}
\theta^{*}=\left(1-\frac{c}{\|Y\|}\right)Y,
\end{equation}
where $c$\ is a positive constant which will be specified below. It will be shown that such an estimator allows one to obtain an explicit upper bound for the quadratic risk in case of the regression model \eqref{sec:int.4} with a conditionally Gaussian noise. Theorem \ref{Le.sec:Gas.1} in Section \ref{sec:Gas} claims that the estimator \eqref{sec:int.5} outperforms the maximum likelihood estimate $\hat{\theta}_{ML}$\ uniformly in $\theta$\ from any compact set $\Theta\subset\mathbb{R}^{p}$\ for any dimension $p$\ starting from two. In Section \ref{sec:Per}, we apply the estimator \eqref{sec:int.5} to solve the problem of improved parametric estimation in the regression model in continuous time with a non-Gaussian noise.

The rest of the paper is organized as follows. In Section \ref{sec:Gas}, we impose some conditions on the random covariance matrix $\mathcal{D}(\mathcal{G})$\ and derive to upper bound for the difference of risks
\begin{equation*}
\Delta(\theta):=R(\theta^{*},\theta)-R(\hat{\theta}_{ML},\theta)
\end{equation*}
corresponding to $\theta^*$\ and $\hat{\theta}_{ML}$\ respectively. In Section \ref{sec:Aut}, the estimate \eqref{sec:int.5} is used for the parameter estimation in a discrete time regression with a Gaussian noise depending on some nuisance parameters. Appendix contains some technical results.


\section{The upper bound for the estimate risk}\label{sec:Gas}

In this section we will derive an upper bound for the risk of estimate \eqref{sec:int.5} under some conditions on the random covariance matrix $\mathcal{D}(\mathcal{G})$.\

Assume that

$(\mathbf{C_1})$\ There exists a positive constant $\lambda_*$,\ such that the minimal eigenvalue of matrix $\mathcal{D}(\mathcal{G})$\ satisfies the inequality
\begin{equation*}
\lambda_{min}(\mathcal{D}(\mathcal{G}))\geq \lambda_*  \quad\mbox{a.s.}
\end{equation*}

$(\mathbf{C_2})$\ The maximal eigenvalue of the matrix $\mathcal{D}(\mathcal{G})$\ is bounded on some compact set $\Theta\subset\mathbb{R}^{p}$\ from above, i.e.
\begin{equation*}
\sup_{\theta\in\Theta}\mathbf{E}_\theta\lambda_{max}(\mathcal{D}(\mathcal{G}))\leq a^*,
\end{equation*}
where $a^*$\ is some known positive constant.

Let denote the difference of the risks of estimate \eqref{sec:int.5} and that of \eqref{sec:int.2} as

\begin{equation*}
\Delta(\theta):=R(\theta^{*},\theta)-R(\hat{\theta}_{ML},\theta).
\end{equation*}

We will need also the following constant
\begin{equation*}
\gamma_p=\dfrac{\sum_{j=0}^{p-2}2^{\frac{j-1}{2}}(-1)^{p-j}\mu^{p-1-j}\Gamma\left(\frac{j+1}{2}\right)-(-\mu)^pI(\mu)}{2^{p/2-1}\Gamma\left(\frac{p}{2}\right)d},
\end{equation*}
where $\mu=d/\sqrt{a^*}$,\
\begin{equation*}
I(a)=\int_0^\infty\frac{\exp(-r^2/2)}{a+r}dr\, \quad\mbox{and} \quad d=\sup\{\|\theta\|: \theta\in\Theta\}.
\end{equation*}

\begin{theorem}\label{Le.sec:Gas.1}

Let the noise $\xi$\ in \eqref{sec:int.4} have a conditionally Gaussian distribution $\mathcal{N}_{p}(0,\mathcal{D}(\mathcal{G}))$\ and its covariance matrix $\mathcal{D}(\mathcal{G})$\ satisfy conditions $(\mathbf{C_1}), (\mathbf{C_2})$\ with some compact set $\Theta\subset\mathbb{R}^{p}$.\ Then the estimator \eqref{sec:int.5} with $c=(p-1)\lambda_{*}\gamma_p$\ dominates the MLE $\hat{\theta}_{ML}$\ for any $p\geq 2$,\ i.e.
\begin{equation*}
\sup_{\theta\in\Theta}\Delta(\theta)\leq-[(p-1)\lambda_{*}\gamma_p]^2.
\end{equation*}

\end{theorem}

\proof First we will establish the lower bound for the random variable $\|Y\|^{-1}$.\

\begin{lemma}\label{Le.sec:Gas.2}
Under the conditions of Theorem 2.1
\begin{equation*}
\inf_{\theta\in\Theta}\mathbf{E}_{\theta}\frac{1}{\|Y\|}\geq\gamma_p.
\end{equation*}
\end{lemma}
The proof of lemma is given in the Appendix.

In order to obtain the upper bound for $\Delta(\theta)$\ we will adjust the argument in the proof of Stein's lemma \cite{St2} to the model \eqref{sec:int.4} with a random covariance matrix.

We consider the risks of MLE and of \eqref{sec:int.5}
\begin{gather*}
R(\hat{\theta}_{ML},\theta)=\mathbf{E}_{\theta}\|\hat{\theta}_{ML}-\theta\|^{2}
=\mathbf{E}_{\theta}(\mathbf{E}\|\hat{\theta}_{ML}-\theta\|^{2}|\mathcal{G})
=\mathbf{E}_{\theta}tr\mathcal{D}(\mathcal{G});\\[2mm]
R(\theta^{*},\theta)=R(\hat{\theta}_{ML},\theta)+\mathbf{E}_{\theta}[\mathbf{E}((g(Y)-1)^{2}\|Y\|^{2}|\mathcal{G})]\\[2mm]
+2\sum_{j=1}^{p}\mathbf{E}_{\theta}[\mathbf{E}((g(Y)-1)Y_{j}(Y_{j}-\theta_j)|\mathcal{G})],
\end{gather*}
where $g(Y)=1-c/\|Y\|$.

Denoting $f(Y)=(g(Y)-1)Y_{j}$\ and applying the conditional density of distribution of a vector $Y$\ with respect to
$\sigma$-algebra $\mathcal{G}$
\begin{equation*}
p_{Y}(x|\mathcal{G})=\frac{1}{(2\pi)^{p/2}\sqrt{\det\mathcal{D}(\mathcal{G})}}
\exp\left(-\frac{(x-\theta)'\mathcal{D}^{-1}(\mathcal{G})(x-\theta)}{2}\right),
 \end{equation*}
one gets
\begin{equation*}
I_{j}:=\mathbf{E}(f(Y)(Y_{j}-\theta_j)|\mathcal{G})=\int_{\mathbb{R}^p}f(x)(x-\theta_j)p_{Y}(x|\mathcal{G})
d x, \quad j=\overline{1,p}.
\end{equation*}

Making the change of variable $u=\mathcal{D}^{-1/2}(\mathcal{G})(x-\theta)$\ and assuming
$\tilde{f}(u)=f(\mathcal{D}^{1/2}(\mathcal{G})u+\theta)$,\ one finds that
\begin{equation*}
I_{j}=\frac{1}{(2\pi)^{p/2}}\sum_{l=1}^{p}\langle
\mathcal{D}^{1/2}(\mathcal{G})\rangle_{jl}\int_{\mathbb{R}^{p}}\tilde{f}(u)u_{l}\exp\left(-\frac{\|u\|^{2}}
{2}\right) d u, \quad j=\overline{1,p},
 \end{equation*}
where $\langle A\rangle_{ij}$\ denotes the $(i,j)$-th\ element of matrix $A$.\ These quantities can be written as
\begin{equation*}
I_{j}=\sum_{l=1}^{p}\sum_{k=1}^{p}\mathbf{E}(<\mathcal{D}^{1/2}(\mathcal{G})>_{jl}
<\mathcal{D}^{1/2}(\mathcal{G})>_{kl}\frac{\partial
f}{\partial u_k}(u)|_{u=Y}|\mathcal{G}), \quad j=\overline{1,p}.
\end{equation*}

Thus, the risk for an estimator \eqref{sec:int.5} takes the form
\begin{gather*}
R(\theta^{*},\theta)=R(\hat{\theta}_{ML},\theta)+\mathbf{E}_{\theta}((g(Y)-1)^{2}\|Y\|^{2})\\
+2\mathbf{E}_{\theta}\left(\sum_{j=1}^{p}\sum_{l=1}^{p}\sum_{k=1}^{p}<\mathcal{D}^{1/2}(\mathcal{G})>_{jl}
<\mathcal{D}^{1/2}(\mathcal{G})>_{kl}\frac{\partial}{\partial
u_{k}}[(g(u)-1)u_j]|_{u=Y}\right).
\end{gather*}

Therefore one has
\begin{equation*}
R(\theta^{*},\theta)=R(\hat{\theta}_{ML},\theta)+\mathbf{E}_{\theta}W(Y),
\end{equation*}
where
\begin{equation*}
W(z)=c^{2}+2c\frac{z'\mathcal{D}(\mathcal{G})z}{\|z\|^3}-2tr\mathcal{D}(\textit{G})c\frac{1}{\|z\|}.
\end{equation*}
This implies that
\begin{equation*}
\Delta(\theta)=\mathbf{E}_{\theta}W(Y).
\end{equation*}

Since $z'Az\leq\lambda_{max}(A)\|z\|^2$,\ one comes to the inequality
\begin{equation*}
\Delta(\theta)\leq c^2-2c\mathbf{E}_{\theta}\frac{tr\mathcal{D}(\mathcal{G})-\lambda_{max}(\mathcal{D}(\mathcal{G}))}{\|Y\|}.
\end{equation*}

From here it follows that
\begin{equation*}
\Delta(\theta)\leq c^2-2c\sum_{i=2}^p\mathbf{E}_{\theta}\frac{\lambda_{i}(\mathcal{D}(\mathcal{G}))}{\|Y\|}.
\end{equation*}
Taking into account the condition $(\mathbf{C_1})$\ and the Lemma ~\ref{Le.sec:Gas.2}, one has
\begin{equation*}
\Delta(\theta)\leq c^2-2(p-1)\lambda_{*}\gamma_p c=:\phi(c).
\end{equation*}

Minimizing the function $\phi(c)$\ with respect to $c$,\ we come to the desired result
\begin{equation*}
\Delta(\theta)\leq -[(p-1)\lambda_{*}\gamma_p]^2.
\end{equation*}
Hence Theorem ~\ref{Le.sec:Gas.1}.

\endproof

\begin{corollary}\label{Le.sec:Gas.3}
Let  in \eqref{sec:int.4} the noise $\xi\sim\mathcal{N}_{p}(0,D)$\ with the positive definite non random covariance matrix $D>0$\ and $\lambda_{min}(D)\geq\lambda_*>0$.\ Then the estimator \eqref{sec:int.5} with $c=(p-1)\lambda_*\gamma_p$\ dominates the MLE for any $p\geq 2$\ and compact set $\Theta\subset\mathbb{R}^{p}$,\ i.e.
\begin{equation*}
\sup_{\theta\in\Theta}\Delta(\theta)\leq-[(p-1)\lambda_{*}\gamma_p]^2.
\end{equation*}

\end{corollary}

\begin{remark}
Note that if $D=\sigma^2I_p$\ then
\begin{equation*}
\sup_{\theta\in\Theta}\Delta(\theta)\leq-[(p-1)\sigma^2\gamma_p]^2.
\end{equation*}
\end{remark}

\begin{corollary}\label{Le.sec:Gas.4}
If $\xi\sim\mathcal{N}_{p}(0,I_p)$ and $\theta=0$\ in model \eqref{sec:int.4} then the risk of estimate \eqref{sec:int.5} is given by the formula
\begin{equation}\label{sec:int.6}
R(0,\theta^*)=p-\left[\frac{(p-1)\Gamma((p-1)/2)}{\sqrt{2}\Gamma(p/2)}\right]^2=:r_p.
\end{equation}
\end{corollary}

By applying the Stirling's formula for the Gamma function
\begin{equation*}
\Gamma(x)=\sqrt{2\pi}x^{x-1/2}\exp(-x)\left(1+o(1)\right)
\end{equation*}
one can check that $r_p\rightarrow 0.5$\ as $p\rightarrow\infty$.\ The behavior of the risk \eqref{sec:int.6} for small values of $p$\ is shown in Fig.1. It will be observed that in this case the risk of the James--Stein estimate remains constant for all $p\geq3$,\ i.e. 
\begin{equation*}
R(0,\hat{\theta}_{JS})=2
\end{equation*}
and the risk of the MLE $\hat{\theta}_{ML}$\ is equal to $p$\ and tends to infinity as $p\rightarrow\infty$.\

\begin{figure}[t]
\centering
\includegraphics[width=0.9\textwidth]{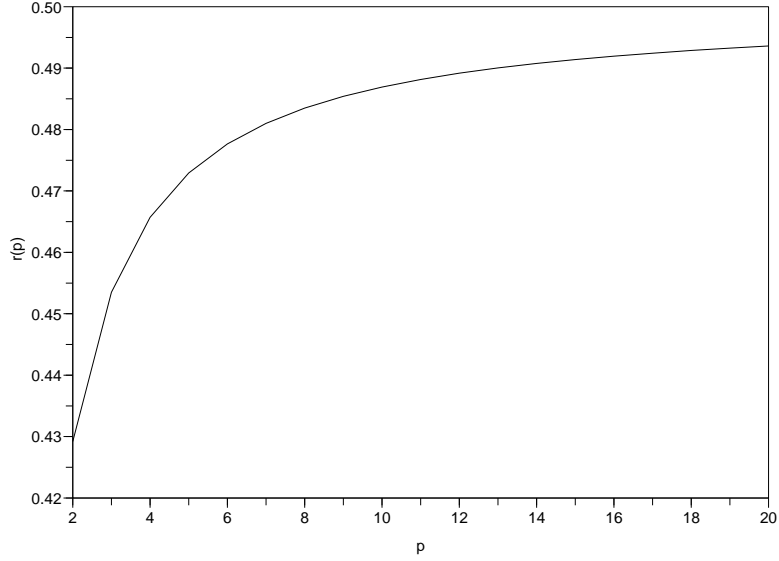}
\caption{Risk of $\theta^*$\ at $\theta=0$.\ }
\label{Risk_figure}
\end{figure}

\section{Improved estimation in a non-Gaussian Ornstein--Uhlenbeck--Levy regression model}\label{sec:Per}

In this section we apply the proposed estimate \eqref{sec:int.5} to a non-Gaussian continuous time regression model. Let observations $(y_t)_{0\leq t\leq n}$\ obey the equation
\begin{equation}\label{sec:Per.1}
d y_{t}=\sum_{j=1}^p\theta_j\phi_j(t)d t+d \xi_{t},\quad 0\leq t\leq n.
 \end{equation}
Here a vector $\theta=(\theta_1,...,\theta_p)'$\ of unknown parameters from some compact set $\Theta\subset\mathbb{R}^{p}$.\
Assume that $(\phi_j(t))_{j\geq1}$\ is a one-periodic $\mathbb{R}_+\rightarrow\mathbb{R}$ functions, bounded and orthonormal  in $L_2[0,1]$.\ The noise $(\xi_{t})_{t\geq 0}$\ in \eqref{sec:Per.1} is a non-Gaussian Ornstein--Uhlenbeck process given by stochastic differential equation
\begin{equation}\label{sec:Per.3}
d \xi_{t}=a \xi_{t}d t+d u_{t},
 \end{equation}
where $a\leq 0$\ is unknown parameter, $(u_{t})_{t \geq 0}$\ is a Levy process satisfying the equation
\begin{equation}\label{sec:Per.4}
u_{t}=\varrho_{1}w_{t}+\varrho_{2}z_{t}.
 \end{equation}
Here $\varrho_{1}, \varrho_{2}$\ are unknown constants, $(w_t)_{t\geq 0}$\ is a standard Brownian motion,
$(z_t)_{t\geq 0}$\ is a compound Poisson process defined as 
\begin{equation}\label{sec:Per.5}
z_{t}=\sum_{j=1}^{N_t} Y_{j},
 \end{equation}
where $(N_t)_{t\geq 0}$\ is a Poisson process with unknown intensity $\lambda>0$ and $(Y_j)_{j\geq1}$\ is a sequence of i.i.d. Gaussian random variables with parameters (0,1).

The problem is to estimate the unknown vector parameter $\theta$\ on the basis of observations $(y_t)_{0\leq t\leq n}$.

Let $\mathcal{G}=\sigma\{N_{t}, t\geq 0\}$\ denote the $\sigma$-algebra\ generated by the Poisson process.

It will be noted that the model \eqref{sec:Per.1} is conditionally Gaussian given the $\sigma$-algebra\ $\mathcal{G}$.\ Therefore one can use estimate \eqref{sec:int.5} to obtain an improved estimate of the unknown vector parameter $\theta$.\ To this end we have to reduce the initial continuous time regression model \eqref{sec:Per.1} to a discrete time model of the form \eqref{sec:int.4}.

The quality of an estimator $\tilde{\theta}$\ will be measured by the quadratic risk
\begin{equation*}
R(\theta,\tilde{\theta})=\mathbf{E}_{\theta}\|\theta-\tilde{\theta}\|^2.
\end{equation*}

A commonly used estimator of an unknown vector $\theta$\ in model \eqref{sec:Per.1} is the least squares estimate (LSE)
$\hat{\theta}=(\hat{\theta}_1,...,\hat{\theta}_p)'$\ with 
\begin{equation*}
\hat{\theta_j}=\frac{1}{n}\int_0^n\phi_j(t)dy_t,\quad j=\overline{1,p}.
\end{equation*}
From here taking into account \eqref{sec:Per.1}, one has
\begin{equation}\label{sec:Per.7}
\hat{\theta}=\theta+\varepsilon_n\zeta(n),
\end{equation}
where $\varepsilon_n=n^{-1/2}$\ and $\zeta(n)$\ is the random vector with coordinates
\begin{equation*}
\zeta_j(n)=\frac{1}{\sqrt{n}}\int_0^n\phi_j(t)d\xi_t.
\end{equation*}
Note that the vector $\zeta(n)$\ has a conditionally Gaussian distribution with a zero mean and conditional covariance matrix $V_{n}(\mathcal{G})=cov(\zeta(n),\zeta(n)'|\mathcal{G})$\ with the elements
\begin{equation*}
v_{ij}(n)=\mathbf{E}(\zeta_{i}(n)\zeta_{j}(n)\mid\mathcal{G}).
\end{equation*}

Thus the initial problem of estimating parameter $\theta$\ in \eqref{sec:Per.1} can be reduced to the that of estimating parameter $\theta$\ in conditionally Gaussian regression model \eqref{sec:Per.7}.

\begin{theorem}\label{Le.sec:Per.1}
Let the regression model be given by the equations \eqref{sec:Per.1}--\eqref{sec:Per.5}, $\varrho_{1}>0$.\ Then, for any $n\geq 2$\ and $p\geq 1$,\ the estimator of $\theta$\
\begin{equation*}
\theta^*=\left(1-\frac{\varrho_{1}^{2}(p-1)\gamma_p}{n\Vert\hat{\theta}\Vert}\right)\hat{\theta},
\end{equation*}
dominates the LSE $\hat{\theta}$:\
\begin{equation*}
\sup_{\theta\in\Theta}\Delta(\theta)\leq -\left[\frac{\varrho_{1}^{2}(p-1)\gamma_p}{n}\right]^2.
\end{equation*}

\end{theorem}

To proved this theorem one can apply Theorem \ref{Le.sec:Gas.1} and it suffices to check conditions $(\mathbf{C_1})$,\ $(\mathbf{C_2})$\ on the matrix $V_{n}(\mathcal{G})$.\ The proof of conditions $(\mathbf{C_1})$\ and $(\mathbf{C_2})$\ is given in the Appendix.

\section{Improved estimation in an autoregression}\label{sec:Aut}

In this section we consider the problem of improved estimating the unknown mean of a multivariate normal distribution when the dispersion matrix is unknown and depends on some nuisance parameters.

Let in \eqref{sec:int.4} the noise $\xi=(\xi_1,...,\xi_d)'$,\ be described by a Gaussian autoregression process
\begin{equation}\label{sec:Aut.1}
\xi_k=a\xi_{k-1}+\varepsilon_k,\ k=\overline{1,p},
\end{equation}
where $|a|<1$,\ $\mathbf{E}\xi_0=0$\ and $\varepsilon_1,...,\varepsilon_p$\ are independent Gaussian (0,1) random variables. Assume that the parameter $a$\ in \eqref{sec:Aut.1} is unknown and belongs to interval $[-\alpha,\alpha]$,\ where $0<\alpha<1$\ is known number.

It is easy to check that the covariance of the noise $\xi$\ has the form

\begin{equation*}
\mathcal{D}(a)=\frac{1}{1-a^2}
\left(
\begin{array}{llll}
1 & a &...& a^{p-1} \\[2mm]
     a & 1 &...& a^{p-2} \\[2mm]
     &  \ddots & &   \\[2mm]
a^{p-1} & a^{p-2} &...& 1
\end{array}
\right)
\end{equation*}

\begin{proposition}
Let $\xi$ in \eqref{sec:int.4} be specified by \eqref{sec:Aut.1} with $a\in [-\alpha,\alpha]$.\ Then for any $p>1/(1-\alpha)^2$\ the MLE is dominated by the estimator
\begin{equation*}
\theta^*=\left(1-\left(p-\frac{1}{(1-\alpha)^2}\right)\frac{\gamma_p}{\|Y\|}\right)Y
\end{equation*}
and
\begin{equation*}
\sup_{\theta\in\Theta}\Delta(\theta)\leq -\left(p-\frac{1}{(1-\alpha)^2}\right)^2\gamma_p^2.
\end{equation*}

\end{proposition}

\proof One has that $tr\mathcal{D}(a)=p/(1-a^2)$.\ Now we find the estimation of the maximal eigenvalue of matrix $\mathcal{D}(a)$.\ By definition
\begin{equation*}
\lambda_{max}(\mathcal{D}(a))=\sup_{\Vert z\Vert=1}z'\mathcal{D}(a)z
\end{equation*}
one has
\begin{gather*}
z'\mathcal{D}(a)z=\sum_{i=1}^{p}\sum_{j=1}^{p}<\mathcal{D}(a)>_{ij}z_iz_j=\frac{1}{1-a^2}\left(1+2\sum_{i=1}^{p-1}\sum_{j=1}^{p-i}a^jz_iz_{j+i}\right)\\
=\frac{1}{1-a^2}\left(1+2\sum_{j=1}^{p-1}a^j\sum_{i=1}^{p-j}z_jz_{i+j}\right).
\end{gather*}
By applying the Cauchy--Bunyakovskii inequality we obtain that
\begin{equation*}
\lambda_{max}(\mathcal{D}(a))\leq\frac{1}{1-\alpha^2}\left(1+2\sum_{j=1}^{\infty}\alpha^j\right)=\frac{1}{(1-\alpha)^2}.
\end{equation*}
Thus,
\begin{equation*}
tr\mathcal{D}(a)-\lambda_{max}(\mathcal{D}(a))\geq p-\frac{1}{(1-\alpha)^2}.
\end{equation*}
Hence, taking into account the Theorem \ref{Le.sec:Gas.1} we come to assertion of Proposition.

\endproof

\section{Conclusions}

In this paper we propose a new type improved estimation procedure. The main difference from the well-known James--Stein estimate is that in the dominator in the corrected term we take the first power of the observation norm $\|Y\|$.\ This allow us to improve estimation with respect to MLE begining with any dimension $p\geq2$.\ Moreover, we apply this procedure to the estimation problem for the non-Gaussian Ornstein--Uhlenbeck--Levy regression model.

\section{Appendix}\label{sec:App}

6.1. Proof of the Lemma \ref{Le.sec:Gas.2}.

\proof From \eqref{sec:int.4} one has
\begin{equation*}
J=\mathbf{E}_{\theta}\frac{1}{\|Y\|}=\mathbf{E}_{\theta}\frac{1}{\|\theta+\xi\|}\geq\mathbf{E}_{\theta}\frac{1}{d+\|\xi\|}.
\end{equation*}
Using a repeated conditional expectation and since the random vector $\xi$\ is distributed conditionally normal with a zero mean, then
\begin{equation*}
J\geq\mathbf{E}_\theta\frac{1}{(2\pi)^{p/2}\sqrt{det\mathcal{D}(\mathcal{G})}}
\int_{\mathbb{R}^p}\frac{\exp(-x'\mathcal{D}(\mathcal{G})^{-1}x/2)}{d+\|x\|}dx.
\end{equation*}
Making the change of variable $u=\mathcal{D}(\mathcal{G})^{-1/2}x$\ and applying the estimation
$u'\mathcal{D}(\mathcal{G})u\leq\lambda_{max}(\mathcal{D}(\mathcal{G}))\|u\|^2$\ we find
\begin{equation*}
J\geq\frac{1}{(2\pi)^{p/2}}
\int_{\mathbb{R}^p}\frac{\exp(-\|u\|^2/2)}{d+\sqrt{\lambda_{max}(\mathcal{D}(\mathcal{G}))}\|u\|}du.
\end{equation*}

Further making the spherical changes of the variables yields

\begin{equation*}
J\geq\frac{1}{2^{p/2-1}\Gamma(p/2)}\mathbf{E}_\theta\int_{0}^\infty\frac{r^{p-1}\exp(-r^2/2)}{d+\sqrt{\lambda_{max}(\mathcal{D}(\mathcal{G}))}r}dr.
\end{equation*}

From here by Jensen and Cauchy--Bunyakovskii inequalities and by the condition $(\mathbf{C_2})$\ we obtain
\begin{equation*}
J\geq\frac{\mu}{2^{p/2-1}\Gamma(p/2)d}
\int_{0}^\infty\frac{r^{p-1}\exp(-r^2/2)}{\mu+r}dr=\gamma_p.
\end{equation*}

This leads to the assertion of Lemma \ref{Le.sec:Gas.2}.

\endproof 

6.2. The proof of conditions $(\mathbf{C_1})$\ and $(\mathbf{C_2})$\ on the matrix $V_{n}(\mathcal{G})$.\

The elements of matrix $V_{n}(\mathcal{G})$\ can be written as \cite{KoPer}
\begin{multline}\label{sec:App.1}
v_{ij}(n)=\frac{\varrho_{1}^{2}}{n}\int_{0}^{n}\phi_i(t)\phi_j(t)dt\\
+\frac{\varrho_{1}^{2}}{2n}\int_{0}^{n}\left(\phi_i(t)\varepsilon_{\phi_j}(t)+\phi_j(t)\varepsilon_{\phi_i}(t)\right)dt
+\frac{\varrho_{2}^{2}}{n}\sum_{l\geq 1}\phi_i(T_l)\phi_j(T_l)\chi_{(T_{l}\leq n)}\\
+\frac{\varrho_{2}^{2}}{n}\sum_{l\geq1}\int_{0}^{n}\left(\phi_i(t)L_{\phi_j}(t-T_{l},T_{l})+\phi_j(t)L_{\phi_i}(t-T_{l},T_{l})\right)\chi_{(T_{l}\leq t)}dt,
\end{multline}
where
\begin{gather*}
\varepsilon_{g}(t)=a\int_{0}^{t}\exp(a(t-s))g(s)(1+\exp(2as))d s,\\
L_{g}(x,y)=a\exp(ax)\left(g(y)+a\int_{0}^{x}\exp(as)g(s+y) d
s\right)
\end{gather*}
and $(T_l)_{l\geq 1}$\ are the jump times of the Poisson process $(N_t)_{t\geq 0}$,\ i.e.
\begin{equation*}
T_{l}=\inf\{t\geq 0: N_{t}=l\}.
\end{equation*}

\begin{lemma}\label{Lem.sec:App.1}
Let $(\xi_{t})_{t\geq 0}$\ be defined by \eqref{sec:Per.3} with $a\leq 0$.\ Then a matrix $V_{n}(\mathcal{G})=(v_{ij}(\phi))_{1\leq i,j\leq p}$\ with elements defined by \eqref{sec:App.1}, satisfy the following inequality a.s.
\begin{equation*}
\inf_{n\geq 1}\inf_{\|z\|=1}z'V_{n}(\mathcal{G})z\geq\varrho_{1}^{2}.
\end{equation*}
\end{lemma}

\proof Notice that by \eqref{sec:App.1} one can the matrix
$V_{n}(\mathcal{G})$\ present as
\begin{equation*}
V_{n}(\mathcal{G})=\varrho_{1}^{2}I_p+F_n+B_{n}(\mathcal{G}),
\end{equation*}
where $F_n$\ is non random matrix with elements
\begin{equation*}
f_{ij}(n)=\frac{\varrho_{1}^{2}}{2n}\int_{0}^{n}\left(\phi_i(t)\varepsilon_{\phi_j}(t)+\phi_j(t)\varepsilon_{\phi_i}(t)\right)dt
\end{equation*}
and $B_{n}(\mathcal{G})$\ is a random matrix with elements
\begin{multline*}
b_{ij}(n)=\varrho_{2}^{2}\sum_{l\geq
1}[\phi_i(T_l)\phi_j(T_l)\chi_{(T_{l}\leq n)}\\
+\int_{0}^{n}\left(\phi_i(t)L_{\phi_j}(t-T_{l},T_{l})+\phi_j(t)L_{\phi_i}(t-T_{l},T_{l})\right)\chi_{(T_{l}\leq
t)}d t ].
 \end{multline*}
This implies that
\begin{equation*}
z'V_{n}(\mathcal{G})z=\varrho_{1}^{2}z'z+z'F_nz+z'B_{n}(\mathcal{G})z\geq \varrho_{1}^{2}z'z,
\end{equation*}
therefore
\begin{equation*}
\inf_{\|z\|=1}z'V_{n}(\mathcal{G})z\geq\varrho_{1}^{2}
\end{equation*}
and we come to the assertion of Lemma \ref{Lem.sec:App.1}.

\endproof

\begin{lemma}\label{Lem.sec:App.2}
Let $(\xi_{t})_{t\geq 0}$\ be defined by \eqref{sec:Per.3} with $a\leq 0$.\ Then a maximal eigenvalue of the matrix $V_{n}(\mathcal{G})=(v_{ij}(n))_{1\leq i,j\leq p}$\ with elements defined by \eqref{sec:App.1}, satisfy the following inequality
\begin{equation*}
\sup_{n\geq1}\sup_{\theta\in\Theta}\mathbf{E}_\theta\lambda_{max}(V_{n}(\mathcal{G}))\leq M p\varrho^*,
\end{equation*}
where $\varrho^*=\varrho_{1}^{2}+\lambda\varrho_{2}^{2}$\ and $M>0$.\

\end{lemma}

\proof One has 
\begin{equation*}
\mathbf{E}_\theta\lambda_{max}(V_{n}(\mathcal{G}))\leq\mathbf{E}_\theta tr(V_{n}(\mathcal{G}))=\sum_{j=1}^p\mathbf{E}_\theta \zeta_j^2(n)=\frac{\varrho^*}{n}\sum_{j=1}^p\tau_{j}(n),
\end{equation*}
where
\begin{multline*}
\tau_{j}(n)=\int_0^n\phi_j^2(t)dt \\
+a\int_0^n\phi_j(t)\int_0^t\exp(a(t-u))\phi_j(u)(1+\exp(2au))dudt.
\end{multline*}
Since the $(\phi_j)_{1\leq j\leq p}$ is a one-periodic orthonormal functions therefore the first integral is equal to $n$\ and in view of the inequality $\max|\phi_j(t)|\leq K$ for any $a\leq 0$\ 
\begin{multline*}
a\int_0^n\phi_j(t)\int_0^t\exp(a(t-u))\phi_j(u)(1+\exp(2au))dudt  \\
\leq 2K^2|a|\int_0^n\int_0^t\exp(a(t-u))dudt\leq 2K^2n.
\end{multline*}
From here denoting $M=1+2K^2$ we obtain that
\begin{equation*}
\mathbf{E}_\theta\lambda_{max}(V_{n}(\mathcal{G}))\leq M p\varrho^*.
\end{equation*}
Hence Lemma \ref{Lem.sec:App.2}.

\endproof

Thus the matrix $V_{n}(\mathcal{G})$\ is positive definite and satisfies for any compact set
$\Theta\subset\mathbb{R}^{p}$,\ the conditions $(\mathbf{C_1})$\ and $(\mathbf{C_2})$\ with $\lambda_*=\varrho_1^2$\ and $a^*=Mp\varrho_*$.\


\medskip

\begin{thebibliography}{100}






\bibitem{BeHa}
J.O. Berger, L.R. Haff, A class of minimax estimators of a normal mean vector for arbitrary quadratic loss and unknown covariance matrix, Statist. Decisions 1 (1983) 105-129.




\bibitem{EfMo}
B. Efron, C. Morris, Families of minimax estimators of the mean of a multivariate normal distribution. Ann. Statist. 4 (1976) 11-21.

\bibitem{Fou}
D. Fourdrinier, Statistique inf\'erentielle, Dunod, 2002, p. 336.

\bibitem{FoPe}
D. Fourdrinier, S. Pergamenshchikov, Improved selection model method for the regression with dependent noise, Ann. of the Inst. of Statist. Math., 59 (3) (2007) 435-464.



\bibitem{FoSt2}
D. Fourdrinier, W.E. Strawderman, E. William, A unified and generalized set of shrinkage bounds on minimax Stein estimates, J. Multivariate Anal. 99 (2008) 2221-2233.


\bibitem{Gl2}
L.J. Gleser, Minimax estimators of a normal mean vector for arbitrary quadratic loss and unknown covariance matrix, Ann. Statist. 14 (1986) 1625-1633.

\bibitem{JS}
W. James, C. Stein, Estimation with quadratic loss, in: Proceedings of the Fourth Berkeley Symposium on Mathematics Statistics and Probability, Vol. 1, University of California Press, Berkeley, 1961, pp. 361-380.

\bibitem{KoPer}
V. Konev, S. Pergamenchtchikov, Efficient robust nonparametric estimation in a semimartingale regression model, http://hal.archives-ouvertes.fr/hal-00526915/fr/ (2010).




\bibitem{St2}
C. Stein, Estimation of the mean of a multivariate normal distribution, Ann. Statist. 9(6) (1981) 1135-1151.


\end{thebibliography}
\end{document}